\input amstex
\documentstyle{amsppt}
\pagewidth{5.4in}
\pageheight{7.6in}
\magnification=1200
\TagsOnRight
\NoRunningHeads
\topmatter
\title
\bf Uniform Sobolev inequalities for manifolds evolving by Ricci flow
\endtitle
\author
Shu-Yu Hsu
\endauthor
\affil
Department of Mathematics\\
National Chung Cheng University\\
168 University Road, Min-Hsiung\\
Chia-Yi 621, Taiwan, R.O.C.\\
e-mail:syhsu\@math.ccu.edu.tw
\endaffil
\date
Aug 7, 2007
\enddate
\address
e-mail address:syhsu\@math.ccu.edu.tw
\endaddress
\abstract
Let $M$ be a compact n-dimensional manifold, $n\ge 2$, with metric 
$g(t)$ evolving by the Ricci flow $\partial g_{ij}/\partial t=-2R_{ij}$ 
in $(0,T)$ for some $T\in\Bbb{R}^+\cup\{\infty\}$ with $g(0)=g_0$.  
Let $\lambda_0(g_0)$ be the first eigenvalue of the operator $-\Delta_{g_0}
+\frac{R(g_0)}{4}$ with respect to $g_0$. We extend a recent result of R.~Ye 
and prove uniform logarithmic Sobolev inequality and uniform Sobolev 
inequalities along the Ricci flow for any $n\ge 2$ when either $T<\infty$ 
or $\lambda_0(g_0)>0$. As a consequence we extend Perelman's local 
$\kappa$-noncollapsing result along the Ricci flow for any $n\ge 2$ in 
terms of upper bound for the scalar curvature when either $T<\infty$ or 
$\lambda_0(g_0)>0$.
\endabstract
\endtopmatter
\NoBlackBoxes
\define \pd#1#2{\frac{\partial #1}{\partial #2}}
\define \1{\partial}
\define \2{\overline}
\define \3{\varepsilon}
\define \4{\widetilde}
\define \5{\underline}
\define \ov#1#2{\overset{#1}\to{#2}}
\define \oa#1{\overset{a}\to{#1}}
\document

Recently there is a lot of studies on Ricci flow on manifolds because it 
is an important tool in understanding the geometry of manifolds 
\cite{H1--3}, \cite{Hs1--3}, \cite{KL}, \cite{MT}, \cite{P1}, 
\cite{P2}. On the other hand given any compact n-dimensional manifold $M$, 
$n\ge 2$, with a fixed metric $g$ it is known that Sobolev inequalities 
hold \cite{He}. More specifically for any $q\in [1,n)$ and $p$ satisfying
$$
\frac{1}{p}=\frac{1}{q}-\frac{1}{n}\tag 1
$$
there exists a minimal constant $C_{p,q}(M,g)>0$ such that the following holds
$$
\biggl (\int_M|u|^pdV_g\biggr )^{\frac{1}{p}}
\le C_{p,q}(M,g)\biggl (\biggl (\int_M|\nabla u|^q\,dV_g\biggr )^{\frac{1}{q}}
+\frac{1}{\text{vol}_g(M)^{\frac{1}{n}}}
\biggl (\int_M|u|^qdV_g\biggr )^{\frac{1}{q}}\biggr )\tag 2
$$ 
for any $u\in W^{1,q}(M,g)$.
Note that since the above inequality is scale invariant with respect to the 
metric $g$, $C_{p,q}(M,g)=C_{p,q}(M,\lambda g)$ for any $\lambda>0$.
Let $M$ be a n-dimensional compact manifold, $n\ge 2$, with a metric $g(t)$ 
evolving by the the Ricci flow
$$
\frac{\1}{\1 t}g_{ij}=-2R_{ij}\quad\text{ in }M\times (0,T)\tag 3
$$ 
for some $T\in\Bbb{R}^+\cup\{\infty\}$ with initial metric $g(0)=g_0$. 
Because of the importance of Sobolev inequalities in the analysis on manifolds,
it is interesting to know whether there exists a constant $C>0$ in place of 
$C_{p,q}(M,g(t))$ such that some variants of (2) hold uniformly along the 
Ricci flow. In \cite{CL} S.C.~Chang and P.~Lu proved that the Yamabe constant
which corresponds to the case $q=2$, $p=2n/(n-2)$, and $n\ge 3$ in (1) has a 
nonnegative derivative at the initial time of the Ricci flow under a technical 
assumption. 

Recently R.~Ye \cite{Ye} by using the monotonicity property of the 
$W$-functional of Perelman \cite{P1} proved the uniform logarithmic Sobolev 
inequality for $n\ge 3$ along the Ricci flow. He also proved that when 
$T<\infty$, then there exist constants $A>0$, $B>0$, such that the uniform 
Sobolev inequality 
$$
\biggl (\int_M|u|^{\frac{2n}{n-2}}dV_{g(t)}\biggr )^{\frac{n-2}{n}}
\le A\int_M\biggl (|\nabla u|^2+\frac{R}{4}u^2\biggr )\,dV_{g(t)}
+B\int_Mu^2\,dV_{g(t)}
$$ 
holds along the Ricci flow for any $n\ge 3$. Similar result was also proved 
by R.~Ye \cite{Ye} when $\lambda_0(g_0)>0$ where $\lambda_0(g_0)$
is the first eigenvalue of the operator $-\Delta_{g_0}+\frac{R(g_0)}{4}$
with respect to $g_0$. As a consequence
R.~Ye \cite{Ye} also obtain the $\kappa$-noncollapsing result in terms of
upper bound on the scalar curvature along the Ricci flow for $n\ge 3$ when 
either $T<\infty$ or $\lambda_0(g_0)>0$.

In this paper we will use a modification of the technique of R.~Ye \cite{Ye},
Qi.S.~Zhang \cite{Z} and E.B.~Davis \cite{D} to prove the uniform logarithmic 
Sobolev inequality and uniform Sobolev inequalities (8) along the Ricci 
flow for any $n\ge 2$, $q\in (1,n)$ and $p$ satisfying (1) when either 
$T<\infty$ or $\lambda_0(g_0)>0$. As a consequence we extend Perelman's 
local $\kappa$-noncollapsing result along the Ricci flow in 
terms of upper bound on the scalar curvature for 
$n\ge 2$ when either $T<\infty$ or $\lambda_0(g_0)>0$. 
Note that Perelman's local $\kappa$-noncollapsing result is stated in terms
of upper bound for $|\text{Rm}|$ \cite{P1} while the $\kappa$-noncollapsing 
result in \cite{Ye} and in this paper is in terms of upper bound on the 
scalar curvature.

We will assume that $(M,g(t))$, $0\le t<T$, is a n-dimensional compact 
manifold, $n\ge 2$, with metric $g(t)$ evolving by the Ricci flow (3)
for the rest of the paper. For any metric $g$ on $M$, let $R(g)$ be 
the scalar curvature of $(M,g)$, $R_-=-\min (0,R)$, and $R(x,t)
=R(g(t))(x)$. For any $0\le t<T$, $r>0$, $x\in M$,
we let $B_t(x,r)$ be the geodesic ball at $x$ with radius $r$ with respect to
the metric $g(t)$ and we let $vol_t(A)$ be the volume of any measurable 
set $A\subset M$ with respect to the metric $g(t)$. 

Note that since 
$$
\frac{\1}{\1 t}\sqrt{\text{det}(g_{ij}(t))}=-R,
$$
$L^p(M,g(t))=L^p(M,g_0)$ and $W^{1,p}(M,g(t))=W^{1,p}(M,g_0)$for all 
$1\le p\le\infty$. Hence we can write $L^p(M)$ and $W^{1,p}(M)$ for 
$L^p(M,g(t))$ and $W^{1,p}(M,g(t))$ for any $t\in [0,T)$ and $1\le p\le\infty$.
For any metric $g$ on $M$, $1\le p<\infty$, $f\in L^p(M,g)$, let
$$
\|f\|_{p,g}=\biggl (\int_M|f|^p\,dV_{g}\biggr )^{\frac{1}{p}}.
$$
When there is no ambiguity we will write $\|f\|_p$ for $\|f\|_{p,g}$.

\proclaim{\bf Theorem 1}
Let $n\ge 2$, $q\in [1,n)$ and $p$ satisfy (1). For any $v\in W^{1,q}(M,g)$, 
$\|v\|_q=1$, we have 
$$
\int_M|v|^q\log v^2\,dV_g\le\frac{2n}{q}\log (C_{p,q}(M,g)(\|\nabla v\|_q
+\text{vol}_g(M)^{-\frac{1}{n}})).
\tag 4
$$
Moreover for any $u\in W^{1,2}(M,g)$, $\|u\|_2=1$,
$$
\int_Mu^2\log u^2\,dV_g\le n\log ((C_{\frac{n\mu}{n-\mu},\mu}(M,g)
((2/\mu)\|\nabla u\|_2+\text{vol}_g(M)^{-\frac{1}{n}}))\tag 5
$$
for any $1\le\mu<2$.
\endproclaim
\demo{Proof}
We will use a modification of the proof of Theorem 3.1 of \cite{Ye} to prove
the lemma. Let $u\in W^{1,q}(M,g)$ satisfy $\|u\|_q=1$. Since $\log x$ is 
concave function, by Jensen's inequality and (2),
$$\align
\int_M|u|^q\log |u|^{p-q}\,dV_g\le&\log\biggl (\int_M|u|^{p-q}\cdot |u|^q
\,dV_g\biggr )
=\log\|u\|_p^p\\
\Rightarrow\quad\int_M|u|^q\log u^2\,dV_g\le&\frac{2p}{p-q}\log (C_{p,q}(M,g)
(\|\nabla u\|_q+\text{vol}_g(M)^{-\frac{1}{n}}\|u\|_q))\\
=&\frac{2n}{q}\log (C_{p,q}(M,g)(\|\nabla u\|_q
+\text{vol}_g(M)^{-\frac{1}{n}})).
\endalign
$$
We now let $\mu\in [1,2)$. For any $u\in W^{1,2}(M,g)$ satisfying $\|u\|_2=1$.
Let $v=|u|^{\frac{2}{\mu}}$. Then $\|v\|_{\mu}=1$. By the Holder inequality,
$$\align
\int_M|\nabla v|^{\mu}\,dV_g
=&\biggl (\frac{2}{\mu}\biggr)^{\mu}
\int_M|u|^{2-{\mu}}|\nabla u|^{\mu}\,dV_g\\
\le&\biggl (\frac{2}{\mu}\biggr)^{\mu}
\biggl (\int_M|\nabla u|^2\,dV_g\biggr )^{\frac{\mu}{2}}
\biggl (\int_Mu^2\,dV_g\biggr )^{1-\frac{\mu}{2}}\\
=&\biggl (\frac{2}{\mu}\biggr)^{\mu}
\biggl (\int_M|\nabla u|^2\,dV_g\biggr )^{\frac{\mu}{2}}\\
\Rightarrow\quad\biggl (\int_M|\nabla v|^{\mu}\,dV_g\biggr )^{\frac{1}{\mu}}
\le&\frac{2}{\mu}\biggl (\int_M|\nabla u|^2\,dV_g\biggr )^{\frac{1}{2}}.
\tag 6
\endalign
$$
Hence $v\in W^{1,\mu}(M,g)$. By putting $q=\mu$, $p=n\mu/(n-\mu)$ and 
$v=|u|^{\frac{2}{\mu}}$ in (4) and using (6) we get (5) and the lemma follows.
\enddemo

By an argument similar to that of Section 3 of \cite{Ye} but with
Theorem 1 replacing Theorem 3.1 of \cite{Ye} in the proof there we get that 
Theorem 3.3 and Theorem 3.5 of \cite{Ye} remain valid for any $n\ge 2$ with 
the constants $C_S(M,g)$ and $vol_g(M)^{-1/n}$ there being replaced by 
$2C_{\frac{n\mu}{n-\mu},\mu}(M,g)/\mu$ and 
$C_{\frac{n\mu}{n-\mu},\mu}(M,g)\text{vol}_g(M)^{-\frac{1}{n}}$ 
respectively for any $1\le\mu<2$. 
Hence Theorem A, Theorem B, Theorem C, and the Corollary on P.1--3 
of \cite{Ye} remain valid for any $n\ge 2$ with the constants $C_S(M,g_0)$ and 
$vol_{g_0}(M)^{-1/n}$ there being replaced by 
$2C_{\frac{n\mu}{n-\mu},\mu}(M,g_0)/\mu$ and 
$C_{\frac{n\mu}{n-\mu},\mu}(M,g_0)\text{vol}_{g_0}(M)^{-\frac{1}{n}}$ 
respectively for any $1\le\mu<2$. In particular we have

\proclaim{\bf Theorem 2}(cf. Theorem A and Theorem C of \cite{Ye})
Let $M$ be a compact n-dimensional manifold, $n\ge 2$, with metric $g(t)$
evolving by the Ricci flow (2) in $(0,T)$ for some 
$T\in\Bbb{R}^+\cup\{\infty\}$ and initial metric $g(0)=g_0$. 
Then there exist constants $C_1>0$, $C_2>0$, depending only on $n$, 
$C_{\frac{n}{n-1},1}(M,g_0)$, the upper bound of $\text{vol}_{g_0}(M)^{-1}$  
and $\max_MR(g_0)_-$ such that 
$$
\int_Mu^2\log u^2\,dV_{g(t)}\le\sigma\int_M\biggl (|\nabla u|+\frac{R}{4}
u^2\biggr )\,dV_{g(t)}-\frac{n}{2}\log\sigma +C_1(t+\frac{\sigma}{4})+C_2
$$
holds for any $\sigma>0$, $0\le t<T$, and $u\in W^{1,2}(M)$ with 
$\|u\|_{2,g(t)}=1$.

If $\lambda_0=\lambda_0(g_0)>0$, then there exists a constant $C_3>0$, 
depending only on $n$, $C_{\frac{n}{n-1},1}(M,g_0)$, a positive lower
bound for $\lambda_0=\lambda_0(g_0)$, the upper bound of 
$\text{vol}_{g_0}(M)^{-1}$  and $\max_MR(g_0)_-$ 
such that 
$$
\int_Mu^2\log u^2\,dV_{g(t)}\le\sigma\int_M\biggl (|\nabla u|+\frac{R}{4}
u^2\biggr )\,dV_{g(t)}-\frac{n}{2}\log\sigma +C_3
$$
holds for any $\sigma>0$, $0\le t<T$, and $u\in W^{1,2}(M)$ with 
$\|u\|_{2,g(t)}=1$.
\endproclaim

For any $0\le s<T$, let 
$$
H_s=-\Delta^s+\frac{(R(g(s))+\max_MR(g(s))_-)}{4}
$$ 
where $\Delta^s$ is the Laplacian of $M$ with respect to the metric $g(s)$. 
Then $H_s$ is a self-adjoint operator on $L^2(M)$. Let 
$$
Q_s(u)=\int_M\biggl (|\nabla u|_{g(s)}^2+\frac{(R(g(s))+\max_MR(g(s))_-)}{4}
u^2\biggr )\,dV_{g(s)}\quad\forall u\in W^{1,2}(M)
$$
be the associated quadratic form of $H_s$. By the maximum principle the 
operator $e^{-H_st}$ is positivity-preserving for all $t\ge 0$ and is a 
contraction on $L^{\infty}(M)$ for all $t\ge 0$. Hence $e^{-H_st}$ is a 
symmetric Markov semigroup (P.22 of \cite{D}). By Theorem 2 and 
the same argument as the proof of Theorem 2.2.7 and Corollary 2.2.8 
of \cite{D} we have the following lemma. 

\proclaim{\bf Lemma 3}
Suppose $n\ge 2$ and either $T<\infty$ or $\lambda_0(g_0)>0$. Then there 
exists a constant $C_4>0$ which depends only on $T$ and the constants $C_1$, 
$C_2$, of Theorem 2 if $T<\infty$ and depends only on the constant $C_3$ of 
Theorem 2 if $\lambda_0(g_0)>0$ such that
$$
\|e^{-H_st}u\|_{\infty}\le C_4t^{-\frac{n}{4}}\|u\|_{2,g(s)}\quad\forall 
0<t\le 1, 0\le s<T, u\in L^2(M).
$$
\endproclaim

By duality we have

\proclaim{\bf Lemma 4}
Suppose $n\ge 2$ and either $T<\infty$ or $\lambda_0(g_0)>0$. Let $C_4>0$ be 
as given in Lemma 3. Then
$$
\|e^{-H_st}u\|_{2,g(s)}\le C_4t^{-\frac{n}{4}}\|u\|_{1,g(s)}\quad\forall 
0<t\le 1, 0\le s<T, u\in L^1(M).
$$
\endproclaim

\proclaim{\bf Corollary 5}
Suppose $n\ge 2$ and either $T<\infty$ or $\lambda_0(g_0)>0$. Let $C_4>0$ be 
as given in Lemma 3 and let $C_5=C_4^2$. Then
$$
\|e^{-H_st}u\|_{\infty}\le C_5t^{-\frac{n}{2}}\|u\|_{1,g(s)}\quad\forall 
0<t\le 1, 0\le s<T, u\in L^1(M).\tag 7
$$
\endproclaim
\demo{Proof}
Since $e^{-H_st}u=e^{-H_s(t/2)}e^{-H_s(t/2)}u$, by Lemma 3 and Lemma 4
for any $0<t\le 1$, $0\le s<T$, $u\in L^1(M)$,
$$
\|e^{-H_st}u\|_{\infty}\le C_4t^{-\frac{n}{4}}\|e^{-H_s(t/2)}u\|_{2,g(s)}
\le C_4^2t^{-\frac{n}{2}}\|u\|_{1,g(s)}.
$$
\enddemo

\proclaim{\bf Theorem 6}
Suppose $n\ge 2$ and either $T<\infty$ or $\lambda_0(g_0)>0$.
Let $1<q<n$ and $p$ satisfy (1). Then there exists a 
constant $A_{p,q}=A_{p,q}(M,g_0)>0$ which depends only on $n$, $T$, 
$C_{\frac{n}{n-1},1}(M,g_0)$, the upper bound of $\text{vol}_{g_0}(M)^{-1}$ 
and $\max_MR(g_0)_-$ if $T<\infty$ and depends only on $n$, 
$C_{\frac{n}{n-1},1}(M,g_0)$, the lower bound of $\lambda_0(g_0)$,
the upper bound of $\text{vol}_{g_0}(M)^{-1}$ and $\max_MR(g_0)_-$ if 
$\lambda_0(g_0)>0$ such that 
$$
\biggl (\int_M|u|^pdV_{g(t)}\biggr )^{\frac{1}{p}}
\le A_{p,q}\biggl (\int_M\biggl (|\nabla u|^2+\frac{R+4+\max_MR(g_0)_-}{4}
u^2\biggr )^{\frac{q}{2}}\,dV_{g(t)}\biggr )^{\frac{1}{q}}\tag 8
$$
holds for any $0\le t<T$ and $u\in W^{1,q}(M)$.
\endproclaim
\demo{Proof}
We will use a modification of the proof of Theorem 2.4.2 of \cite{D} 
(cf. \cite{Ye}, \cite{Z}) to prove
the theorem. Let $\4{H}_s=H_s+1$. Then $\4{H}_s$ is also a symmetric 
Markov semigroup. By Corollary 5,
$$
\|e^{-\4{H}_st}u\|_{\infty}\le C_5t^{-\frac{n}{2}}\|u\|_{1,g(s)}\quad\forall 
0<t\le 1, 0\le s<T, u\in L^1(M)\tag 9
$$
for some constant $C_5>0$ given by (7).
Since $e^{-H_st}$ is a contraction on $L^{\infty}(M)$ for all $t\ge 0$
and $e^{-\4{H}_st}u=e^{1-t}\cdot e^{-H_s(t-1)}e^{-\4{H}_s}u$, by (9),
$$
\|e^{-\4{H}_st}u\|_{\infty}\le e^{1-t}\|e^{-\4{H}_s}u\|_{\infty}
\le C_5e^{1-t}\|u\|_{1,g(s)}\quad\forall t\ge 1, 0\le s<T, u\in L^1(M).
\tag 10
$$
By (9) and (10) there exists a constant $C_6>0$ such that
$$
\|e^{-\4{H}_st}u\|_{\infty}\le C_6t^{-\frac{n}{2}}\|u\|_{1,g(s)}\quad\forall 
t>0, 0\le s<T, u\in L^1(M).\tag 11
$$
By (11) and the Riesz-Thorin interpolation theorem \cite{D},
$$
\|e^{-\4{H}_st}u\|_{\infty}\le C_6^{\frac{1}{q}}t^{-\frac{n}{2q}}
\|u\|_{q,g(s)}\quad\forall t>0, 0\le s<T, u\in L^q(M).\tag 12
$$
By (12) and an argument similar to 
the proof of Theorem 2.4.2 of \cite{D} $\4{H}_s^{-\frac{1}{2}}$ is of weak 
type $(p,q)$ for any $1<q<n$ and $p$ satisfying (1). We now choose $q_1,q_2$
such that $1<q_1<q<q_2<n$ and let $p_i$ satisfy (1) with $p$, $q$ being 
replaced by $p_i$ and $q_i$, $i=1,2$ respectively. Then $p_1<p<p_2$. By the 
above argument 
$\4{H}_s^{-\frac{1}{2}}$ is of weak type $(p_1,q_1)$ and weak type $(p_2,q_2)$.
Hence by the Marcinkiewicz interpolation theorem \cite{S} there exists 
a constant $C>0$ such that
$$
\|\4{H}_s^{-\frac{1}{2}}u\|_{p,g(s)}\le C\|u\|_{q,g(s)}\quad\forall 0\le s<T.
\tag 13
$$ 
Since $R$ satisfy \cite{H1},
$$
R_t=\Delta R+2|\text{Ric}|^2\quad\text{ in }M\times (0,T).
$$ 
Hence $\min_MR(x,t)=-\max_MR_-(x,t)$ is a monotone increasing function of 
$t\in [0,T)$. Then $\max_MR_-(x,t)\le\max_MR(g_0)_-$ for all $0\le t<T$. 
Thus by (13) we get (8) and the theorem follows.
\enddemo

\proclaim{\bf Theorem 7}($\kappa$-noncollapsing result)
Suppose $n\ge 2$ and either $T<\infty$ or $\lambda_0(g_0)>0$. Let $\rho>0$. 
Then for any $x_0\in M$, $0\le t<T$, $0<r\le\rho$, satisfying
$$
R(x,t)\le r^{-2}\quad\forall x\in B_t(x_0,r),\tag 14
$$
the following holds
$$
\text{vol}_t(B_t(x_0,r))\ge\kappa r^n\tag 15
$$
where
$$
\kappa=\min\biggl (\frac{1}{(2^{(n+3q)/q}A_{p,q})^n},
\frac{1}{(\sqrt{2}A_{p,q}(1+(4+\max_MR(g_0)_-)\rho^2)^{1/2})^{n}}
\biggr )\tag 16
$$
for any $1<q<n$ and $p$ satisfying (1) where $A_{p,q}=A_{p,q}(M,g_0)$ 
is given by Theorem 6.
\endproclaim
\demo{Proof}
We will use a modification of the proof of Lemma 2.2 of \cite{He} and 
Lemma 6.1 of \cite{Ye} to prove the theorem. Suppose there exists $x_0\in M$, 
$0\le t<T$, $0<r\le\rho$ satisfying (14) such that (15) does not hold. Then
$$
\text{vol}_t(B_t(x_0,r))<\kappa r^n.
$$
Let $\2{g}=r^{-2}g(t)$. Then
$$
\text{vol}_{\2{g}}(B_{\2{g}}(x_0,1))<\kappa\tag 17
$$ 
and
$$
\2{R}\le 1\quad\text{ in }B_{\2{g}}(x_0,1)\tag 18
$$
where $\2{R}$ is the scalar curvature of $\2{g}$. Let $q\in (1,n)$ and 
$p$ satisfy (1). By Theorem 6 there exists a constant $A_{p,q}=A_{p,q}
(M,g_0)>0$ such that (8) holds. By (8),
$$
\biggl (\int_M|u|^pdV_{\2{g}}\biggr )^{\frac{1}{p}}
\le A_{p,q}\biggl (\int_M\biggl (|\nabla u|_{\2{g}}^2
+\frac{\2{R}+(4+\max_MR(g_0)_-)\rho^2}{4}
u^2\biggr )^{\frac{q}{2}}\,dV_{\2{g}}\biggr )^{\frac{1}{q}}\tag 19
$$
holds for any $u\in W^{1,q}(M,\2{g})$. Then for any $u\in W^{1,q}(M,\2{g})$ 
with compact support in $B_{\2{g}}(x_0,1)$ by (16), (17), (18), (19), and 
the Holder inequality,
$$\align
&\biggl (\int_{B_{\2{g}}(x_0,1)}|u|^pdV_{\2{g}}\biggr )^{\frac{1}{p}}\\
\le&A_{p,q}\biggl (\int_{B_{\2{g}}(x_0,1)}\biggl (|\nabla u|_{\2{g}}^2
+\frac{1+(4+\max_MR(g_0)_-)\rho^2}{4}
u^2\biggr )^{\frac{q}{2}}\,dV_{\2{g}}\biggr )^{\frac{1}{q}}\\
\le&\sqrt{2}A_{p,q}\biggl (\int_{B_{\2{g}}(x_0,1)}
|\nabla u|_{\2{g}}^q\,dV_{\2{g}}\biggr )^{\frac{1}{q}}\\
&\qquad +\sqrt{2}A_{p,q}\frac{(1+(4+\max_MR(g_0)_-)\rho^2)^{\frac{1}{2}}}{2}
\biggl (\int_{B_{\2{g}}(x_0,1)}|u|^q\,dV_{\2{g}}\biggr )^{\frac{1}{q}}\\
\le&\sqrt{2}A_{p,q}\biggl (\int_{B_{\2{g}}(x_0,1)}
|\nabla u|_{\2{g}}^q\,dV_{\2{g}}\biggr )^{\frac{1}{q}}\\
&\qquad +\sqrt{2}A_{p,q}\frac{(1+(4+\max_MR(g_0)_-)\rho^2)^{\frac{1}{2}}}{2}
\kappa^{\frac{1}{n}}\biggl (\int_{B_{\2{g}}(x_0,1)}|u|^p\,dV_{\2{g}}\biggr )^{\frac{1}{p}}\\
\le&\sqrt{2}A_{p,q}\biggl (\int_{B_{\2{g}}(x_0,1)}
|\nabla u|_{\2{g}}^q\,dV_{\2{g}}\biggr )^{\frac{1}{q}}+\frac{1}{2}
\biggl (\int_{B_{\2{g}}(x_0,1)}|u|^p\,dV_{\2{g}}\biggr )^{\frac{1}{p}}.
\endalign
$$
Hence
$$
\biggl (\int_{B_{\2{g}}(x_0,1)}|u|^p\,dV_{\2{g}}\biggr )^{\frac{1}{p}}
\le 2\sqrt{2}A_{p,q}\biggl (\int_{B_{\2{g}}(x_0,1)}
|\nabla u|_{\2{g}}^q\,dV_{\2{g}}\biggr )^{\frac{1}{q}}.\tag 20
$$
For any $0<r_1\le 1$, let $u(x)=r_1-\text{dist}_{\2{g}}(x,x_0)$. By 
the Holder inequality,
$$
\frac{r_1}{2}\text{vol}_{\2{g}}(B_{\2{g}}(x_0,r_1/2))^{\frac{1}{q}}
\le\biggl (\int_{B_{\2{g}}(x_0,r_1)}|u|^q\,dV_{\2{g}}\biggr )^{\frac{1}{q}}
\le\text{vol}_{\2{g}}(B_{\2{g}}(x_0,r_1))^{\frac{1}{n}}
\biggl (\int_{B_{\2{g}}(x_0,1)}|u|^p\,dV_{\2{g}}\biggr )^{\frac{1}{p}}.
\tag 21
$$
Hence by (20) and (21),
$$\align
&\frac{r_1}{2}\text{vol}_{\2{g}}(B_{\2{g}}(x_0,r_1/2))^{\frac{1}{q}}
\le 2\sqrt{2}A_{p,q}\text{vol}_{\2{g}}(B_{\2{g}}(x_0,r_1))^{\frac{n+q}{nq}}\\
\Rightarrow\quad&\text{vol}_{\2{g}}(B_{\2{g}}(x_0,r_1))\ge\biggl (
\frac{r_1}{4\sqrt{2}A_{p,q}}\biggr )^{\frac{nq}{n+q}}
\text{vol}_{\2{g}}(B_{\2{g}}(x_0,r_1/2))^{\frac{n}{n+q}}\quad\forall
0<r_1\le 1.\tag 22
\endalign
$$
By (17), (22), and an argument similar to the proof of Lemma 2.2 of \cite{He}
and \cite{C},
$$\align
&\text{vol}_{\2{g}}(B_{\2{g}}(x_0,r_1))\ge\biggl (\frac{1}{2^{(n+(5q/2))/q}
A_{p,q}}\biggr )^nr_1^n\quad\forall 0<r_1\le 1\\
\Rightarrow\quad&\kappa>\text{vol}_{\2{g}}(B_{\2{g}}(x_0,1))
\ge\biggl (\frac{1}{2^{(n+(5q/2))/q}A_{p,q}}\biggr )^n.
\endalign
$$
This contradicts (16). Hence (15) holds and the theorem follows.
\enddemo

\enddocument